\documentclass[12pt,twoside]{article}
\usepackage{mathrsfs}
\usepackage{amssymb,amsmath,amsthm,epsfig,latexsym}
\usepackage[numbers,sort&compress]{natbib}
\usepackage{txfonts}
\usepackage{amsfonts}
\usepackage{geometry}
\usepackage[latin1]{inputenc}

\newtheorem{theorem}{Theorem}[section]

\newtheorem{example}{Example}

\theoremstyle{definition}
\newtheorem{definition}{Definition}

\numberwithin{equation}{section}

\begin{document}

\renewcommand{\baselinestretch}{1.15}

\begin{center}
{\Large \textbf{Existence and multiplicity of positive weak solutions for a new class of (p; q)-Laplacian systems%
}} \vskip.2in {\small \qquad Rafik Guefaifia$^{2}$ \qquad Jiabin Zuo$^{1,3,4}
$ \qquad Salah Boulaaras\footnote{%
Corresponding author: saleh\_boulaares@yahoo.fr(S.Boulaaras). Co-authors:
zuojiabin88@163.com (J.B. Zuo); rafikik982@gmail.com(R.Guefaifia);
goyal.praveen2011@gmail.com(P. Agarwal).}$^{5,6}$ \qquad Praveen Agarwal$^{7}
$}\\[1mm]
{\scriptsize $^{1}$Faculty of Applied Sciences, Jilin Engineering Normal
University, Changchun 130052, P. R. China}\\[0pt]
{\scriptsize $^{2}$Department of Mathematics, Faculty of Exact Sciences,
University Tebessa, Tebessa, Algeria}\\[0pt]
{\scriptsize $^{3}$College of Science, Hohai University, Nanjing 210098, P.
R. China}\\[0pt]
{\scriptsize $^{4}$Departamento de Matem\'atica, Universidade Estadual de
Campinas, IMECC, Campinas, SP CEP 13083-859, Brazil}\\[0pt]
{\scriptsize $^{5}$Department of Mathematics, College of Sciences and  Arts,
Al-Rass, Qassim University, Saudi Arabia}\\[0pt]
{\scriptsize $^{6}$Laboratory of Fundamental and Applied Mathematics of Oran
(LMFAO), University of Oran 1, Oran, Algeria}\\[0pt]
{\scriptsize $^{7}$Anand International College of Engineering, Near Kanota,
Agra Road, Jaipur-303012, Rajasthan, India}\\[0pt]
\end{center}

\vspace{0.05in}

{\footnotesize \noindent \textbf{Abstract}~ The paper is concerned with the
existence of positive weak solutions for a new class of $\left( p,q\right) $%
-Laplacian elliptic systems in a bounded domain by means of the method of
sub-super solutions. Particularly, we do not need any sign conditions for $%
\gamma \left( 0\right), g\left( 0\right), f\left( 0\right) $ and $h\left(
0\right) $. Moreover, a multiplicity result is obtained when $\gamma \left(
0\right)=g\left( 0\right)=f\left( 0\right)=h\left( 0\right)=0.$ Finally, we
give some examples to verify our main results.\newline
}

{\footnotesize \noindent\textbf{Keywords:~} New elliptic systems: Existence;
Positive solutions; Multiplicity; Sub-super solutions. \vskip 3mm \noindent%
 }

{\footnotesize \vskip.2in }

\section{ Introduction}

{\footnotesize In this paper, we deal with the existence and multiplicity of
positive weak solutions for the following $\left( p,q\right)$-Laplacian
systems
\begin{equation}
\bigskip \left\{
\begin{array}{l}
-\triangle _{p}u-\left\vert u\right\vert ^{p-2}u=\lambda _{1}a\left(
x\right) f\left( v\right) +\mu _{1}\alpha \left( x\right) h\left( u\right)
\text{ in }\Omega , \\
\\
-\triangle _{q}v-\left\vert v\right\vert ^{q-2}v=\lambda _{2}b\left(
x\right) g\left( u\right) +\mu _{2}\beta \left( x\right) \gamma \left(
v\right)\ \ \text{ in }\Omega , \\
\\
u=v=0\qquad \qquad \qquad\qquad \qquad \qquad\qquad \ \ \text{ on }\partial
\Omega ,%
\end{array}%
\right.  \label{1.1}
\end{equation}%
where $\triangle _{s}z=\text{div}\left( \left\vert \nabla z\right\vert
^{s-2}\nabla z\right) ,s>1,\Omega \subset \mathbb{R}^{N}\ \left( N\geq
3\right)$ is a bounded domain with smooth boundary $\partial \Omega ,a\left(
x\right) ,b\left( x\right) ,$ $\alpha \left( x\right) ,\beta \left( x\right)
\in C(\overline{\Omega}) ,$ $\lambda _{1},\lambda _{2},\mu _{1},$ ,$\mu _{2}
\geq0$, $1<p, q<\infty.$ }

{\footnotesize The study of $\left( p,q\right) $-Laplacian systems is a new
and interesting topic. It arises from electrorheological fluids, nonlinear
elasticity theory, etc. (see \cite{A99}-\cite{R10}, \cite{2} \cite%
{5}, \cite{6}, \cite{23} and \cite{27}). A lot of existence results have
been obtained on this class of problems, we refer to ( \cite{X}, \cite{3},
\cite{5}, \cite{8}, \cite{10}, \cite{Y}, \cite{Z}, \cite{17},\cite{20}, \cite{W}, \cite{25}). These problems originate from physical models
and are widely used in many fields such as combustion, mathematical biology,
chemical reactions and so on. Our method is mainly focused on the method of
sub-super solutions (see \cite{17}, \cite{20} for a more detailed discussion). }

{\footnotesize As far as we know, there are very few contributions devoted
to the $\left( p,q\right)$-Laplacian nonlinear elliptic system. Therefore,
with the help of the method of sub-super solutions method, we are inspired
by the paper of \cite{23} in which a new $\left( p,q\right)$-Laplacian
system was discussed and extended our previous results to problem \eqref{1.1}
without assuming any sign conditions for $h\left( 0\right)$, $g\left(
0\right)$, $f\left( 0\right)$, and $\gamma \left( 0\right) $. Furthermore,
when $h\left( 0\right) = f\left( 0\right)= g\left( 0\right)=\gamma \left(
0\right) =0, $ a multiplicity result is given. }

{\footnotesize The outline of the paper is organized as follows: Sec. 2
introduces some definitions and make appropriate assumptions, which will be
used in the body of the paper. In addition, we show the proof of two
important results. Sec. 3, we will illustrate our main results with some
interesting examples. }

\section{ Main results}

{\footnotesize First, in order to get our main results, we will consider the
following hypothesis:\newline
$\left( H1\right) $ There exist $a\left( x\right),$ $\alpha \left( x\right),$
$b\left( x\right),$ $\beta \left( x\right) \in C( \overline{\Omega }) $ such
that
\begin{equation*}
a\left( x\right) \geq a_{1}>0,b\left( x\right) \geq b_{1}>0,\alpha \left(
x\right) \geq \alpha _{1}>0,\beta \left( x\right) \geq \beta _{1}>0.
\end{equation*}
}

{\footnotesize \noindent$\left( H2\right) $ Let $f,$ $g,$ $h,$ $\gamma \in
C^{1}([0,\infty ))$ be monotone functions satisfying
\begin{equation*}
\lim_{s\rightarrow +\infty }f\left( s\right) =\lim_{s\rightarrow +\infty
}g\left( s\right) =\lim_{s\rightarrow +\infty }h\left( s\right)
=\lim_{s\rightarrow +\infty }\gamma \left( s\right) =+\infty .
\end{equation*}
\noindent$\left( H3\right) $ $\lim_{_{s\rightarrow +\infty }}\frac{f\left(
M\left( g\left( s\right) \right) ^{\frac{1}{q-1}}\right) }{s^{p-1}}=0,$ \ $%
\forall M>0.$\newline
}

{\footnotesize \noindent$\left( H4\right) $ $\lim_{_{s\rightarrow +\infty }}%
\frac{h\left( s\right) }{s^{p-1}}=\lim_{_{s\rightarrow +\infty }}\frac{%
\gamma \left( s\right) }{s^{q-1}}=0.$\newline
}

{\footnotesize Next, we define weak solutions and sub-super solutions in $%
\left( p,q\right)$-Laplacian elliptic systems. }

\begin{definition}
{\footnotesize Let $\left( u,v\right) \in W^{1,p}\left( \Omega \right) \cap
C( \overline{\Omega }) \times W^{1,q}\left( \Omega \right) \cap C( \overline{%
\Omega }) ,$ we say that $\left( u,v\right) $ is a weak solution of problem (%
\ref{1.1}), if }
\end{definition}

{\footnotesize
\begin{equation*}
\begin{array}{c}
\int\limits_{\Omega }\left\vert \nabla u\right\vert ^{p-2}\nabla u.\nabla
\xi dx-\int\limits_{\Omega }\left\vert u\right\vert ^{p-2}u.\xi dx=\lambda
_{1}\int\limits_{\Omega }a\left( x\right) f\left( v\right) \xi dx+\mu
_{1}\int\limits_{\Omega }\alpha \left( x\right) h\left( u\right) \xi dx\text{
in }\Omega , \\
\\
\int\limits_{\Omega }\left\vert \nabla v\right\vert ^{q-2}\nabla v.\nabla
\zeta dx-\int\limits_{\Omega }\left\vert v\right\vert ^{q-2}v.\zeta
dx=\lambda _{2}\int\limits_{\Omega }b\left( x\right) g\left( u\right) \zeta
dx+\mu _{2}\int\limits_{\Omega }\beta \left( x\right) \gamma \left( v\right)
\zeta dx\text{ in }\Omega%
\end{array}%
\end{equation*}
}

{\footnotesize for all $\left( \xi ,\zeta \right) \in W_{0}^{1,p}\left(
\Omega \right) \times W_{0}^{1,q}\left( \Omega \right) .$ }

\begin{definition}
{\footnotesize The nonnegative functions $( \underline{u},\underline{v}) ,(
\overline{u},\overline{v}) $ in $W^{1,p}\left( \Omega \right) \cap C(%
\overline{\Omega}) \times W^{1,q}\left( \Omega \right) \cap C(\overline{%
\Omega}) $ are called a weak subsolution and supersolution of problem (\ref%
{1.1}) if \ they satisfy $( \underline{u},\underline{v}) ,\left( \overline{u}%
,\overline{v}\right) =$ $\left( 0,0\right) $ on $\partial \Omega $ }
\end{definition}

{\footnotesize
\begin{equation*}
\begin{array}{c}
\int\limits_{\Omega }\left\vert \nabla \underline{u}\right\vert ^{p-2}\nabla
\underline{u}.\nabla \xi dx-\int\limits_{\Omega }\left\vert \underline{u}%
\right\vert ^{p-2}\underline{u}.\xi dx\leq \lambda _{1}\int\limits_{\Omega
}a\left( x\right) f\left( \underline{v}\right) \xi dx+\mu
_{1}\int\limits_{\Omega }\alpha \left( x\right) h\left( \underline{u}\right)
\xi dx\text{ in }\Omega , \\
\\
\int\limits_{\Omega }\left\vert \nabla \underline{v}\right\vert ^{q-2}\nabla
\underline{v}.\nabla \zeta dx-\int\limits_{\Omega }\left\vert \underline{v}%
\right\vert ^{q-2}\underline{v}.\zeta dx\leq \lambda _{2}\int\limits_{\Omega
}b\left( x\right) g\left( \underline{u}\right) \zeta dx+\mu
_{2}\int\limits_{\Omega }\beta \left( x\right) \gamma \left( \underline{v}%
\right) \zeta dx\text{ in }\Omega%
\end{array}%
\end{equation*}
}

{\footnotesize and }

{\footnotesize
\begin{equation*}
\begin{array}{c}
\int\limits_{\Omega }\left\vert \nabla \overline{u}\right\vert ^{p-2}\nabla
\overline{u}.\nabla \xi dx-\int\limits_{\Omega }\left\vert \overline{u}%
\right\vert ^{p-2}\overline{u}.\xi dx\geq \lambda _{1}\int\limits_{\Omega
}a\left( x\right) f\left( \overline{v}\right) \xi dx+\mu
_{1}\int\limits_{\Omega }\alpha \left( x\right) h\left( \overline{u}\right)
\xi dx\text{ in }\Omega , \\
\\
\int\limits_{\Omega }\left\vert \nabla \overline{v}\right\vert ^{q-2}\nabla
\overline{v}.\nabla \zeta dx-\int\limits_{\Omega }\left\vert \overline{v}%
\right\vert ^{q-2}\overline{v}.\zeta dx\geq \lambda _{2}\int\limits_{\Omega
}b\left( x\right) g\left( \overline{u}\right) \zeta dx+\mu
_{2}\int\limits_{\Omega }\beta \left( x\right) \gamma \left( \overline{v}%
\right) \zeta dx\text{ in }\Omega%
\end{array}%
\end{equation*}
}

{\footnotesize for any $\left( \xi ,\zeta \right) \in W_{0}^{1,p}\left(
\Omega \right) \times W_{0}^{1,q}\left( \Omega \right) .$\newline
}

{\footnotesize In what follows, we shall establish the following the
existence result. }

\begin{theorem}
{\footnotesize Let $\left( H1\right) -\left( H4\right) $ hold. If $\lambda
_{1}+\mu _{1}$ and $\lambda _{2}+\mu _{2}$ are big enough, then problem (\ref%
{1.1}) processes a positive weak solution. }
\end{theorem}

{\footnotesize
\begin{proof}
We will show that there exist a positive weak subsolution$%
( \underline{u},\underline{v}) \in W^{1,p}\left( \Omega \right)
\cap C(\overline{\Omega}) \times W^{1,q}\left( \Omega \right)
\cap C(\overline{\Omega}) $ and a supersolution $(\overline{u},\overline{v}) \in W^{1,p}\left( \Omega \right) \cap
C(\overline{\Omega}) \times W^{1,q}\left( \Omega \right) \cap
C(\overline{\Omega}) $ of $  \left( 1.1\right) $ such that $%
\underline{u}\leq \overline{u},$ $\underline{v}$ $\leq \overline{v}.$ Moreover, $( \underline{u},\underline{v}) ,( \overline{u},%
\overline{v}) $ satisfy $(\underline{u},\underline{v})
=\left( 0,0\right) =\left( \overline{u},\overline{v}\right) $ on $\partial
\Omega .$

Let  $\sigma _{r}$ be the first eigenvalue of $-\triangle _{r}$,  and $\phi _{r}>0$ the corresponding eigenfunction with $\left\Vert \phi _{r}\right\Vert =1$ for $%
r=p,q.$ There exist $m,\eta ,\delta >0$ such that $\left\vert \nabla \phi _{r}\right\vert ^{r}-\sigma _{r}$ $\phi _{r}\geq m$
on $\overline{\Omega _{\delta }}=\left\{ x\in \Omega ,d\left( x,\partial
\Omega \right) \leq \delta \right\} $ and $\phi _{r}\geq \eta $ on $\Omega
\backslash \overline{\Omega }_{\delta }$ for $r=p,q.$ Taking $k_{0}>0$ such
that $a_{1}f\left( t\right) ,$ $\alpha _{1}h\left( t\right) ,$ $b_{1}g\left(
t\right) ,$ $\beta _{1}\gamma \left( t\right) >-k_{0}.$

First, we claim that

\begin{equation}
(\underline{u},\underline{v}):=\left(\left[ \frac{\left( \lambda _{1}+\mu _{1}\right) k_{0}}{m}%
\right] ^{1\diagup p-1}\left( \frac{p-1}{p}\right) \phi _{p}^{p\diagup p-1}, \left[ \frac{\left( \lambda _{2}+\mu _{2}\right) k_{0}}{m}%
\right] ^{1\diagup q-1}\left( \frac{q-1}{q}\right) \phi _{q}^{q\diagup q-1}\right)
\tag{2.1}
\end{equation}%

is a subsolution of problem $(1.1)$ when $\lambda _{1}+\mu _{1}$ and $\lambda
_{2}+\mu _{2}$ are big enough. Taking the test function $\xi \left( x\right) \in
W_{0}^{1,p}\left( \Omega \right) $ with $\xi \left( x\right) \geq 0.$ Thus,
from $\left( H1\right) $ we get

\begin{eqnarray*}
\int\limits_{\Omega }\left\vert \nabla \underline{u}\right\vert ^{p-2}\nabla
\underline{u}.\nabla \xi dx-\int\limits_{\Omega }\left\vert \underline{u}%
\right\vert ^{p-2}\underline{u}.\xi dx &\leq &\int\limits_{\Omega
}\left\vert \nabla \underline{u}\right\vert ^{p-2}\nabla \underline{u}%
.\nabla \xi dx \\
&=&\left( \frac{\left( \lambda _{1}+\mu _{1}\right) k_{0}}{m}\right)
\int\limits_{\Omega }\left\{ \sigma _{p}\phi _{p}^{p}-\left\vert \nabla \phi
_{p}\right\vert ^{p}\right\} \xi dx \\
&=&\left( \frac{\left( \lambda _{1}+\mu _{1}\right) k_{0}}{m}\right)
\int\limits_{\Omega _{\delta }}\left\{ \sigma _{p}\phi _{p}^{p}-\left\vert
\nabla \phi _{p}\right\vert ^{p}\right\} \xi dx \\
&&+\left( \frac{\left( \lambda _{1}+\mu _{1}\right) k_{0}}{m}\right)
\int\limits_{\Omega \backslash \overline{\Omega }_{\delta }}\left\{ \sigma
_{p}\phi _{p}^{p}-\left\vert \nabla \phi _{p}\right\vert ^{p}\right\} \xi dx.
\end{eqnarray*}

We have known that  $\left\vert \nabla \phi
_{r}\right\vert ^{r}-\sigma _{r}$ $\phi _{r}\geq m$ for $s=p,q$, on $\overline{\Omega _{\delta }}$. Also on $%
\Omega \backslash \overline{\Omega }_{\delta }$ $\phi _{r}\geq \eta $ for $%
r=p,q.$ If $\lambda _{1}+\mu _{1}$ and $\lambda _{2}+\mu _{2}$ are big enough in
the definition of $\underline{u}$, $\underline{v}$, then by $\left( H2\right) $ we get
\begin{equation}
a_{1}f\left( \underline{v}\right) ,\alpha _{1}h\left( \underline{u}\right)
,b_{1}g\left( \underline{u}\right) ,\beta _{1}\gamma \left( \underline{v}%
\right) \geq \frac{k_{0}}{m}\max \left\{ \sigma _{p},\sigma _{q}\right\}.
\tag{2.3}
\end{equation}

Therefore,
\begin{eqnarray*}
\int\limits_{\Omega }\left\vert \nabla \underline{u}\right\vert ^{p-2}\nabla
\underline{u}.\nabla \xi dx-\int\limits_{\Omega }\left\vert \underline{u}%
\right\vert ^{p-2}\underline{u}.\xi dx &\leq &\left( \frac{\left( \lambda
_{1}+\mu _{1}\right) k_{0}}{m}\right) \int\limits_{\Omega _{\delta }}\left\{
\sigma _{p}\phi _{p}^{p}-\left\vert \nabla \phi _{p}\right\vert ^{p}\right\}
\xi dx \\
&&+\left( \frac{\left( \lambda _{1}+\mu _{1}\right) k_{0}}{m}\right)
\int\limits_{\Omega \backslash \overline{\Omega }_{\delta }}\left\{ \sigma
_{p}\phi _{p}^{p}-\left\vert \nabla \phi _{p}\right\vert ^{p}\right\} \xi dx
\\
&\leq &-\left( \lambda _{1}+\mu _{1}\right) k_{0}\int\limits_{\Omega
_{\delta }}\xi dx+\left( \frac{\left( \lambda _{1}+\mu _{1}\right) k_{0}}{m}%
\right) \int\limits_{\Omega \backslash \overline{\Omega }_{\delta }}\sigma
_{p}\xi dx \\
&\leq &\int\limits_{\Omega _{\delta }}\left[ \lambda _{1}a\left( x\right)
f\left( \underline{v}\right) +\mu _{1}\alpha \left( x\right) h\left(
\underline{u}\right) \text{ }\right] \xi dx \\
&&+\int\limits_{\Omega \backslash \overline{\Omega }_{\delta }}\left[
\lambda _{1}a\left( x\right) f\left( \underline{v}\right) +\mu _{1}\alpha
\left( x\right) h\left( \underline{u}\right) \text{ }\right] \xi dx \\
&=&\int\limits_{\Omega }\left[ \lambda _{1}a\left( x\right) f\left(
\underline{v}\right) +\mu _{1}\alpha \left( x\right) h\left( \underline{u}%
\right) \text{ }\right] \xi dx.
\end{eqnarray*}

Similarly,%
\begin{equation*}
\int\limits_{\Omega }\left\vert \nabla v\right\vert ^{q-2}\nabla v.\nabla
\zeta dx-\int\limits_{\Omega }\left\vert v\right\vert ^{q-2}v.\zeta dx\leq
\int\limits_{\Omega }\left[ \lambda _{2}b\left( x\right) g\left( \underline{u%
}\right) +\mu _{2}\beta \left( x\right) \gamma \left( \underline{v}\right) %
\right] \zeta dx.
\end{equation*}

Thus $( \underline{u},\underline{v}) $ is a subsolution of
problem $\left( 1.1\right) .$

Next, let $\omega _{r}$ be a unique
positive solution of
\begin{equation*}
\left\{
\begin{array}{c}
-\triangle _{r}\omega _{r}=1\text{ in }\Omega , \\
\\
\omega _{r}=0\ \ \ \ \ \ \ \text{ on }\partial \Omega .%
\end{array}%
\right.
\end{equation*}

for $r=p,q.$ We denote

\begin{equation}
\overline{u}=\frac{C}{\nu _{p}}\left( \frac{\lambda _{1}\left\Vert
a\right\Vert _{\infty }+\mu _{1}\left\Vert \alpha \right\Vert _{\infty }}{%
1-\nu _{p}^{p-1}}\right) ^{\frac{1}{p-1}}\omega _{p},  \tag{2.4}
\end{equation}

\begin{equation}
\overline{v}=\left[ \left( \frac{\lambda _{2}\left\Vert b\right\Vert
_{\infty }+\mu _{2}\left\Vert \beta \right\Vert _{\infty }}{1-\nu _{q}^{q-1}}%
\right) g\left( C\left( \frac{\lambda _{1}\left\Vert a\right\Vert _{\infty
}+\mu _{1}\left\Vert \alpha \right\Vert _{\infty }}{1-\nu _{p}^{p-1}}\right)
^{\frac{1}{p-1}}\right) ^{\frac{1}{q-1}}\right] \omega _{q},  \tag{2.5}
\end{equation}

where $\nu _{r}=\left\Vert \omega _{r}\right\Vert _{\infty },$ $r=p,q$ and $%
C>0$ is big enough. We claim that $\left(
\overline{u},\overline{v}\right) $ is a supersolution of $(1.1)$ such that $%
\left( \overline{u},\overline{v}\right) \geq ( \underline{u},\underline{%
v}) .$

According to $(H3)-(H4)$, we can make $C$ big enough so that
\begin{equation}
\begin{array}{l}
\left( \frac{C}{\nu _{p}}\right) ^{p-1}\geq f\left( \left[ \left( \frac{%
\lambda _{2}\left\Vert b\right\Vert _{\infty }+\mu _{2}\left\Vert \beta
\right\Vert _{\infty }}{1-\nu _{q}^{q-1}}\right) g\left( C\left( \frac{%
\lambda _{1}\left\Vert a\right\Vert _{\infty }+\mu _{1}\left\Vert \alpha
\right\Vert _{\infty }}{1-\nu _{p}^{p-1}}\right) ^{\frac{1}{p-1}}\right) ^{%
\frac{1}{q-1}}\right] \omega _{q}\right)
+\mu _{1}h\left( \frac{\lambda _{1}\left\Vert a\right\Vert _{\infty }+\mu
_{1}\left\Vert \alpha \right\Vert _{\infty }}{1-\nu _{p}^{p-1}}\right) ^{%
\frac{1}{p-1}}\omega _{p}.%
\end{array}
\tag{2.6}
\end{equation}

Hence
\begin{equation*}
\int\limits_{\Omega }\left\vert \nabla \overline{u}\right\vert ^{p-2}\nabla
\overline{u}.\nabla \xi dx-\int\limits_{\Omega }\left\vert \overline{u}%
\right\vert ^{p-2}\overline{u}.\xi dx=\left( \frac{C}{\nu _{p}}\right)
^{p-1}\left( \lambda _{1}\left\Vert a\right\Vert _{\infty }+\mu
_{1}\left\Vert \alpha \right\Vert _{\infty }\right) \int\limits_{\Omega }\xi
dx.
\end{equation*}

Using $\left( 2.6\right) $
\begin{equation}
\begin{array}{l}
\int\limits_{\Omega }\left\vert \nabla \overline{u}\right\vert ^{p-2}\nabla
\overline{u}.\nabla \xi dx-\int\limits_{\Omega }\left\vert \overline{u}%
\right\vert ^{p-2}\overline{u}.\xi dx \\
\geq \lambda _{1}\left\Vert a\right\Vert _{\infty }f\left( \left[ \left(
\frac{\lambda _{2}\left\Vert b\right\Vert _{\infty }+\mu _{2}\left\Vert
\beta \right\Vert _{\infty }}{1-\nu _{q}^{q-1}}\right) g\left( C\left( \frac{%
\lambda _{1}\left\Vert a\right\Vert _{\infty }+\mu _{1}\left\Vert \alpha
\right\Vert _{\infty }}{1-\nu _{p}^{p-1}}\right) ^{\frac{1}{p-1}}\right) ^{%
\frac{1}{q-1}}\right] \omega _{q}\right) \int\limits_{\Omega }\xi dx
+\mu _{1}\left\Vert \alpha \right\Vert _{\infty }\int\limits_{\Omega
}h\left( C\left( \frac{\lambda _{1}\left\Vert a\right\Vert _{\infty }+\mu
_{1}\left\Vert \alpha \right\Vert _{\infty }}{1-\nu _{p}^{p-1}}\right) ^{%
\frac{1}{p-1}}\right) \xi dx \\
\geq \int\limits_{\Omega }\left[ \lambda _{1}a\left( x\right) f\left(
\overline{v}\right) +\mu _{1}\alpha \left( x\right) h\left( \overline{u}%
\right) \right] \text{ }\xi dx.%
\end{array}
\tag{2.7}
\end{equation}

Next

\begin{equation}
\begin{array}{l}
\int\limits_{\Omega }\left\vert \nabla \overline{v}\right\vert ^{q-2}\nabla
\overline{v}.\nabla \zeta dx-\int\limits_{\Omega }\left\vert \overline{v}%
\right\vert ^{q-2}\overline{v}.\zeta dx \\
=\left\{ \left( \lambda _{2}\left\Vert b\right\Vert _{\infty }+\mu
_{2}\left\Vert \beta \right\Vert _{\infty }\right) g\left( C\left( \frac{%
\lambda _{1}\left\Vert a\right\Vert _{\infty }+\mu _{1}\left\Vert \alpha
\right\Vert _{\infty }}{1-\nu _{p}^{p-1}}\right) ^{\frac{1}{p-1}}\right)
\right\} \omega _{q}\int\limits_{\Omega }\xi dx \\
\geq \left[ \lambda _{2}\left\Vert b\right\Vert _{\infty }g\left( C\left(
\frac{\lambda _{1}\left\Vert a\right\Vert _{\infty }+\mu _{1}\left\Vert
\alpha \right\Vert _{\infty }}{1-\nu _{p}^{p-1}}\right) ^{\frac{1}{p-1}%
}\right) +\mu _{2}\left\Vert \beta \right\Vert _{\infty }g\left( C\left(
\frac{\lambda _{1}\left\Vert a\right\Vert _{\infty }+\mu _{1}\left\Vert
\alpha \right\Vert _{\infty }}{1-\nu _{p}^{p-1}}\right) ^{\frac{1}{p-1}%
}\right) \right] \int\limits_{\Omega }\xi dx.%
\end{array}
\tag{2.8}
\end{equation}

According to $(H4)$ and choose $C$ big enough, we obtain
\begin{equation*}
g\left( C\left( \frac{\lambda _{1}\left\Vert a\right\Vert _{\infty }+\mu
_{1}\left\Vert \alpha \right\Vert _{\infty }}{1-\nu _{p}^{p-1}}\right) ^{%
\frac{1}{p-1}}\right) \geq \gamma \left( \left[ \left( \frac{\lambda
_{2}\left\Vert b\right\Vert _{\infty }+\mu _{2}\left\Vert \beta \right\Vert
_{\infty }}{1-\nu _{q}^{q-1}}\right) g\left( C\left( \frac{\lambda
_{1}\left\Vert a\right\Vert _{\infty }+\mu _{1}\left\Vert \alpha \right\Vert
_{\infty }}{1-\nu _{p}^{p-1}}\right) ^{\frac{1}{p-1}}\right) ^{\frac{1}{q-1}}%
\right] \left\Vert \omega _{q}\right\Vert _{\infty }\right).
\end{equation*}

Then from $(2.7)$ we get
\begin{equation}
\begin{array}{l}
\int\limits_{\Omega }\left\vert \nabla \overline{v}\right\vert ^{q-2}\nabla
\overline{v}.\nabla \zeta dx-\int\limits_{\Omega }\left\vert \overline{v}%
\right\vert ^{q-2}\overline{v}.\zeta dx \\
\geq \lambda _{2}\left\Vert b\right\Vert _{\infty }g\left( C\left( \frac{%
\lambda _{1}\left\Vert a\right\Vert _{\infty }+\mu _{1}\left\Vert \alpha
\right\Vert _{\infty }}{1-\nu _{p}^{p-1}}\right) ^{\frac{1}{p-1}}\right)
+\mu _{2}\left\Vert \beta \right\Vert _{\infty }\gamma \left( \left\{ \left(
\frac{\lambda _{2}\left\Vert b\right\Vert _{\infty }+\mu _{2}\left\Vert
\beta \right\Vert _{\infty }}{1-\nu _{q}^{q-1}}\right) g\left( C\left( \frac{%
\lambda _{1}\left\Vert a\right\Vert _{\infty }+\mu _{1}\left\Vert \alpha
\right\Vert _{\infty }}{1-\nu _{p}^{p-1}}\right) ^{\frac{1}{p-1}}\right)
\right\} ^{\frac{1}{q-1}}\left\Vert \omega _{q}\right\Vert _{\infty }\right)
\\
\geq \int\limits_{\Omega }\left[ b\left( x\right) g\left( \overline{u}%
\right) +\mu _{2}\beta \left( x\right) \gamma \left( \overline{v}\right) %
\right] \zeta dx.%
\end{array}
\tag{2.9}
\end{equation}

According to $(2.7)$ and $(2.8)$, we can conclude that $\left( \overline{u},%
\overline{v}\right) $ is a supersolution of $(1.1)$. Further $\overline{u}%
\geq \underline{u}$ and $\overline{v}\geq \underline{v}$ \ for $C$ \ big enough.
Thus, we get  a solution $\left( u,v\right) \in W^{1,p}\left( \Omega
\right) \cap C(\overline{\Omega}) \times W^{1,q}\left( \Omega
\right) \cap C(\overline{\Omega}) $ of $\left( 1.1\right) $
with \ $\underline{u}\leq u\leq \overline{u}$, and $\underline{v}\leq v\leq
\overline{v}.$ The proof of Theorem $1.1$ is complete.
\end{proof}
}

{\footnotesize Now we show that the more general system $(1.1)$ possesses at
least two distinct positive solutions. }

\begin{theorem}
{\footnotesize Suppose that the conditions $\left( H1\right) -\left(
H4\right) $ hold. Let $f,g,h,$ and $\gamma $ the function be smooth enough
around zero with $f\left( 0\right) =h\left( 0\right) =g\left( 0\right)
=\gamma \left( 0\right) =0=f^{\left( k\right) }\left( 0\right) =h^{\left(
k\right) }\left( 0\right) =g^{\left( l\right) }\left( 0\right) =\gamma
^{\left( l\right) }\left( 0\right) $ for $k=1,2,....\left[ p-1\right] ,$ $%
l=1,2,....\left[ q-1\right] , $ where $\left[ s\right] $ denotes the integer
part of $s$. Then, problem $\left( 1.1\right) $ processes at least two
positive solutions when $\lambda _{i}+\mu _{i}$ are big enough; $i=1,2.$ }
\end{theorem}

{\footnotesize
\begin{proof}
For problem $\left( 1.1\right) $, we will look for a strict supersolution $\left( \zeta _{1},\zeta _{2}\right), $ a subsolution $\left( \psi _{1},\psi_{2}\right),$  a supersolution $\left( z_{1},z_{2}\right) $, and a strict subsolution $\left( \omega _{1},\omega _{2}\right) ,$  such
that $\left( \psi _{1},\psi _{2}\right) \leq \left( \zeta _{1},\zeta
_{2}\right) \leq \left( z_{1},z_{2}\right) ,$ $\left( \psi _{1},\psi
_{2}\right) \leq \left( \omega _{1},\omega _{2}\right) \leq \left(
z_{1},z_{2}\right) ,$ and $\left( \omega _{1},\omega _{2}\right) \nleqslant
\left( \zeta _{1},\zeta _{2}\right) .$ Then, problem $\left( 1.1\right) $ processes  at
least three distinct solutions $\left( u_{i},v_{i}\right) ,$ $i=1,2,3,$ such
that
$$\left( u_{1},v_{1}\right) \in \left[ \left( \psi _{1},\psi _{2}\right)
,\left( \zeta _{1},\zeta _{2}\right) \right] , ~~\left( u_{2},v_{2}\right)
\in \left[ \left( \omega _{1},\omega _{2}\right) ,\left( z_{1},z_{2}\right) %
\right] $$
and
$$\left( u_{3},v_{3}\right) \in \left[ \left( \psi _{1},\psi
_{2}\right) ,\left( z_{1},z_{2}\right) \right] \diagdown \left( \left[
\left( \psi _{1},\psi _{2}\right) ,\left( \zeta _{1},\zeta _{2}\right) %
\right] \cup \left( \omega _{1},\omega _{2}\right) ,\left(
z_{1},z_{2}\right) \right) .$$

It is obvious that $\left( \psi _{1},\psi _{2}\right) =\left( 0,0\right) $
is a (sub)solution. Moreover, we always can  find a big
supersolution $\left( z_{1},z_{2}\right) =\left( \overline{u},\overline{v}%
\right) .$ Consider

\begin{equation}
\left\{
\begin{array}{l}
-\triangle _{p}\omega _{1}-\left\vert \omega _{1}\right\vert ^{p-2}\omega
_{1}=\lambda _{1}a\left( x\right) \widetilde{f}\left( \omega _{2}\right)
+\mu _{1}\alpha \left( x\right) \widetilde{h}\left( \omega _{1}\right) \text{
in }\Omega , \\
\\
-\triangle _{q}\omega _{2}-\left\vert \omega _{2}\right\vert ^{q-2}\omega
_{2}=\lambda _{2}b\left( x\right) \widetilde{g}\left( \omega _{1}\right)
+\mu _{2}\beta \left( x\right) \widetilde{\gamma }\left( \omega _{2}\right)
\text{ in }\Omega , \\
\\
\omega _{1}=\omega _{2}=0\text{ on }\partial \Omega ,%
\end{array}%
\right.  \tag{2.10}
\end{equation}

where  $\widetilde{g}\left( s\right)
=g\left( s\right) -1,$ $\widetilde{\gamma }\left( s\right) =\gamma \left(
s\right) -1,$  $\widetilde{h}%
\left( s\right) =h\left( s\right) -1,$ $\widetilde{f}\left( s\right) =f\left( s\right) -1.$  Then by Theorem 1.1, when $\lambda _{i}+\mu
_{i} $ are big enough, we know that the problem $\left( 2.10\right) $ processes a
solution $\left( \omega _{1},\omega _{2}\right)>0 $  $i=1,2.$ It is clear that $\left( \omega _{1},\omega
_{2}\right) $ is a strict subsolution of problem $\left( 1.1\right) $.

In the end, we will find a strict supersolution $\left( \zeta _{1},\zeta _{2}\right) .$

Let $\phi _{p},$ $\phi _{q}$  be the corresponding eigenfunction with respect to operators $%
\triangle _{p}$ and $\triangle _{q}$  and  there exist  $C_{1}>0$ and $C_{2}>0$ such that
\begin{equation}
\phi _{p}\leq C_{1}\phi _{q}\text{ \ \ and \ \ }\phi _{q}\leq C_{2}\phi _{p}.
\tag{2.11}
\end{equation}

Let $\left( \zeta _{1},\zeta _{2}\right) =\left( \rho \phi _{p},\rho \phi
_{q}\right), $ where $\rho >0,$
\begin{eqnarray*}
G_{p}\left( x\right):=\left( \sigma _{p}-1\right) x^{p-1}-\lambda
_{1}f\left( C_{2}x\right) -\mu _{1}h\left( x\right) ~~\text{and} ~~G_{q}\left( x\right):=\left( \sigma _{q}-1\right) x^{q-1}-\lambda
_{2}g\left( C_{1}x\right) -\mu _{2}\gamma \left( x\right) .
\end{eqnarray*}

Note that $G_{p}\left( 0\right) =G_{q}\left( 0\right) =0,$ $G_{p}^{\left(
k\right) }\left( 0\right) =G_{q}^{\left( l\right) }\left( 0\right) =0$ for $%
k=1,2,....\left[ p-1\right] $ and $l=1,2,....\left[ q-1\right] .$

\begin{equation*}
\left\{
\begin{array}{c}
G_{p}^{\left( p-1\right) }\left( 0\right) >0 ~~\text{and}~~ G_{q}^{\left( q-1\right)
}\left( 0\right) >0 ~~~~~~~~~~~~~~~~\text{if} ~~~p,q \in \mathbb{Z^{+}}, \\
\\
\lim_{r\rightarrow
+\infty }G_{p}^{\left( \left[ p\right] \right) }\left( r\right) =+\infty
=\lim_{r\rightarrow +\infty }G_{p}^{\left( \left[ q\right] \right) }\left(
r\right)\ \ ~~\text{if} ~~~p,q \notin \mathbb{Z^{+}}.%
\end{array}%
\right.
\end{equation*}

Hence, there exists $\theta$ such that $G_{q}\left( x\right) >0$ and $G_{p}\left( x\right) >0$ for any $x\in
(0,\theta ].$ So, for $0<\rho \leq \theta $ we get

\begin{equation*}
\left( \sigma _{p}-1\right) \zeta _{1}^{p-1}=\left( \sigma _{p}-1\right)
\left( \rho \phi _{p}\right) ^{p-1}>\lambda _{1}f\left( C_{2}\rho \phi
_{p}\right) -\mu _{1}h\left( \rho \phi _{p}\right).
\end{equation*}

By $\left( 2.11\right) $ and the monotonicity of function $f$, we obtain \
\begin{equation}
\begin{array}{l}
\left( \sigma _{p}-1\right) \zeta _{1}^{p-1}=\left( \sigma _{p}-1\right)
\left( \rho \phi _{p}\right) ^{p-1}>\lambda _{1}f\left( C_{2}\rho \phi
_{p}\right) -\mu _{1}h\left( \rho \phi _{p}\right)
\geq \lambda _{1}f\left( \rho \phi _{q}\right) -\mu _{1}h\left( \rho \phi
_{p}\right)
=\lambda _{1}f\left( \zeta _{2}\right) -\mu _{1}h\left( \zeta _{1}\right)%
\text{ \ \ }
\end{array}
\tag{2.12}
\end{equation}
for any $x\in \Omega$. In the same way, we also have

\begin{equation}
\begin{array}{l}
\left( \sigma _{q}-1\right) \zeta _{2}^{q-1}=\left( \sigma _{q}-1\right)
\left( \rho \phi _{q}\right) ^{q-1}>\lambda _{2}g\left( C_{1}\rho \phi
_{q}\right) -\mu _{2}\gamma \left( \rho \phi _{q}\right)
\geq \lambda _{2}g\left( \rho \phi _{p}\right) -\mu _{2}\gamma \left( \rho
\phi _{q}\right)
=\lambda _{2}g\left( \zeta _{1}\right) -\mu _{2}\gamma \left( \zeta
_{2}\right) ,\text{ \ \ }%
\end{array}
\tag{2.13}
\end{equation}
for any $x\in \Omega$.
Making use of  $\left( 2.12\right) $ and $\left( 2.13\right) $, we
obtain
\begin{eqnarray*}
&&\int\limits_{\Omega }\left\vert \nabla \zeta _{1}\right\vert ^{p-2}\nabla
\zeta _{1}.\nabla \xi dx-\int\limits_{\Omega }\left\vert \zeta
_{1}\right\vert ^{p-2}\zeta _{1}.\xi dx \\
&=&\rho ^{p-1}\left\{ \int\limits_{\Omega }\left\vert \nabla \phi
_{p}\right\vert ^{p-2}\nabla \phi _{p}.\nabla \xi dx-\int\limits_{\Omega
}\left\vert \phi _{p}\right\vert ^{p-2}\phi _{p}.\xi dx\right\} \\
&=&\int\limits_{\Omega }\left\{ \sigma _{p}\left( \rho \phi _{p}\right)
^{^{p-1}}-\left( \rho \phi _{p}\right) ^{^{^{p-1}}}\right\} \xi dx,\text{ \
\ \ because }\phi _{p}>0, \\
&=&\int\limits_{\Omega }\left\{ \left( \sigma _{p}-1\right) \left( \rho \phi
_{p}\right) ^{^{p-1}}\right\} \xi dx \\
&=&\int\limits_{\Omega }\left( \sigma _{p}-1\right) \zeta _{1}^{p-1}\xi dx \\
&>&\lambda _{1}\int\limits_{\Omega }f\left( \zeta _{2}\right) \xi dx-\mu
_{1}\int\limits_{\Omega }h\left( \zeta _{1}\right) .\xi dx,
\end{eqnarray*}

Similarly we also have
\begin{eqnarray*}
&&\int\limits_{\Omega }\left\vert \nabla \zeta _{2}\right\vert ^{q-2}\nabla
\zeta _{2}.\nabla \xi dx-\int\limits_{\Omega }\left\vert \zeta
_{2}\right\vert ^{q-2}\zeta _{2}.\xi dx
>\lambda _{2}\int\limits_{\Omega }g\left( \zeta _{1}\right) \xi dx-\mu
_{2}\int\limits_{\Omega }\gamma \left( \zeta _{2}\right) \xi dx.
\end{eqnarray*}

It follows that $\left( \zeta _{1},\zeta _{2}\right) $ is a strict supersolution.
Let $\rho $ small enough so that $\left( \omega _{1},\omega _{2}\right)
\nleqslant \left( \zeta _{1},\zeta _{2}\right) $. So, we can find
solutions $$\left( u_{1},v_{1}\right) \in \left[ \left( \psi _{1},\psi
_{2}\right) ,\left( \zeta _{1},\zeta _{2}\right) \right] ,\left(
u_{2},v_{2}\right) \in \left[ \left( \omega _{1},\omega _{2}\right) ,\left(
z_{1},z_{2}\right) \right]$$
and
$$\left( u_{3},v_{3}\right) \in \left[
\left( \psi _{1},\psi _{2}\right) ,\left( z_{1},z_{2}\right) \right] %
\diagdown \left( \left[ \left( \psi _{1},\psi _{2}\right) ,\left( \zeta
_{1},\zeta _{2}\right) \right] \cup \left( \omega _{1},\omega _{2}\right)
,\left( z_{1},z_{2}\right) \right) .$$

This fact that $\left( u_{1},v_{1}\right) \equiv \left( \psi
_{1},\psi _{2}\right) \equiv \left( 0,0\right) $ can happen due to $\left( \psi _{1},\psi _{2}\right) \equiv \left( 0,0\right) $ is a
solution. So, anyway we can find two
positive solutions $\left( u_{2},v_{2}\right) $ and $\left(
u_{3},v_{3}\right) $. Therefore, we conclude the proof of Theorem 2.2.
\end{proof}
}

\section{Examples}

\begin{example}
{\footnotesize Let%
\begin{eqnarray*}
f\left( x\right) &=&\sum_{i=1}^{m}a_{i}x^{p_{i}}-c_{1},\text{ }g\left(
x\right) =\sum_{j=1}^{n}b_{j}x^{qj}-c_{2} \\
h\left( x\right) &=&\sum_{k=1}^{s}\alpha _{k}x^{r_{k}}-c_{3},\text{ }\gamma
\left( x\right) =\sum_{l=1}^{\tau }\beta _{l}x^{d_{l}}-c_{4},
\end{eqnarray*}
}
\end{example}

{\footnotesize where $d_{j}<\left( q-1\right),$ $r_{k}<\left( p-1\right) ,$ $%
p_{i}q_{j}<\left( p-1\right) \left( q-1\right),$ $a_{i},b_{j},\alpha
_{k},\beta _{l},p_{i},q_{j},r_{k},d_{j},c_{1},c_{2},c_{3},c_{4}\geq 0,$ So,
it is clear that $f,g,h$ and $\gamma $ fulfill the assumptions of Theorem
1.1. }

\begin{example}
{\footnotesize Let
\begin{eqnarray*}
f\left( x\right) &=&\left\{
\begin{array}{c}
x^{p_{1}},\text{ \ \ \ \ \ \ \ \ \ \ \ \ \ \ \ \ \ \ \ \ \ \ \ \ \ \ }x\leq
1, \\
\frac{p_{1}}{p_{2}}x^{p_{2}}+\left( 1-\frac{p_{1}}{p_{2}}\right) ,\text{ \ \
}x>1,%
\end{array}%
\right. ,\text{ }h\left( x\right) =\left\{
\begin{array}{c}
x^{p_{3}},\text{ \ \ \ \ \ \ \ \ \ \ \ \ \ \ \ \ \ \ \ \ \ \ \ \ \ \ }x\leq
1, \\
\frac{p_{_{3}}}{p_{4}}x^{p_{4}}+\left( 1-\frac{p_{_{3}}}{p_{4}}\right) ,%
\text{ \ \ }x>1,%
\end{array}%
\right. \text{\ \ ,} \\
g\left( x\right) &=&\left\{
\begin{array}{c}
x^{q_{1}},\text{ \ \ \ \ \ \ \ \ \ \ \ \ \ \ \ \ \ \ \ \ \ \ \ \ \ \ }x\leq
1, \\
\frac{q_{1}}{q_{2}}x^{q_{2}}+\left( 1-\frac{q_{1}}{q_{2}}\right) ,\text{ \ \
}x>1,%
\end{array}%
\right. ,\text{ }\gamma \left( x\right) =\left\{
\begin{array}{c}
x^{q_{3}},\text{ \ \ \ \ \ \ \ \ \ \ \ \ \ \ \ \ \ \ \ \ \ \ \ \ \ \ }x\leq
1, \\
\frac{q_{_{3}}}{q_{4}}x^{q_{4}}+\left( 1-\frac{q_{_{3}}}{q_{4}}\right) ,%
\text{ \ \ }x>1,%
\end{array}%
\right. \text{\ }
\end{eqnarray*}
}
\end{example}

{\footnotesize where we suppose that
\begin{equation*}
\left\{
\begin{array}{c}
p_{1},p_{3}>p-1 ~~~~~~~~~~~~~~~~\text{if} ~~p\in \mathbb{Z^{+}} , \\
p_{1},p_{3}>\left[ p\right]\ \ ~~~~~~~~~~~~~~~~~~~~~~\text{if} ~~p\notin
\mathbb{Z^{+}}, \\
q_{1},q_{3}>q-1~~~~~~~~~~~~~~~\text{if} ~~q\in \mathbb{Z^{+}}, \\
q_{1},q_{3}>\left[ q\right]~~~~~~~~~~~~~~~~~~~~~~~~~\text{if}~~~~q\notin
\mathbb{Z^{+}}, \\
\end{array}%
\right.
\end{equation*}
$p_{4}<p-1$, $p_{2}q_{2}<\left(p-1\right) \left( q-1\right)$ and $q_{4}<q-1$%
. Clearly, $f,g,h$ and $\gamma $ fulfill all the assumptions of Theorem 2.2.
}

\section*{{\protect\footnotesize {\protect\small Acknowledgement }}}

{\footnotesize The second author was supported by the Fundamental Research
Funds for Central Universities (2019B44914) and the National Key Research
and Development Program of China (2018YFC1508100), the China Scholarship
Council (201906710004). }

\section*{{\protect\footnotesize {\protect\small Availability of data and
materials}}}

{\footnotesize Not applicable. }

\section*{{\protect\footnotesize {\protect\small Competing interests}}}

{\footnotesize The authors declare that there is no conflict of interests
regarding the publication of this manuscript. The authors declare that they
have no competing interests. }

\section*{{\protect\footnotesize {\protect\small Authors' contributions}}}

{\footnotesize The authors contributed equally in this article. They have
all read and approved the final manuscript. }

\clearpage

\topmargin=-1.2in 
\textheight=11in
\oddsidemargin=-10pt
\textwidth=6.5in
\let\raggedleft\raggedright
\vspace{-12pt} \renewcommand{\baselinestretch}{1.1} 

\end{document}